\documentclass{conm-p-l}%[12pt]{amsart}
\usepackage{amscd,amssymb}
\usepackage[mathscr]{eucal}
\begin{document}

\newcommand{\R}{{\mathbb R}}
\newcommand{\C}{{\mathbb C}}
\newcommand{\bR}{\bar{\mathbb R}}
\newcommand{\Z}{{\mathbb Z}}
\newcommand{\B}{{\mathbb B}}
\newcommand{\Id}{{\mathbb I}}
\newcommand{\E}{{\mathbb E}}
\newcommand{\Ee}{{\mathcal E}}
\renewcommand{\L}{{\mathbb L}}
\newcommand{\rht}{\o{H}^2_\R}
\newcommand{\gG}{\Gamma}
\newcommand{\ga}{\alpha}
\newcommand{\gl}{\lambda}
\newcommand{\na}{\ga'}

\newcommand{\Nn}{\mathcal N}
\newcommand{\Ss}{\mathcal S}
\newcommand{\ML}{{\mathcal{ML}(M)}}
\newcommand{\PL}{{\mathcal{PL}(M)}}
\newcommand{\Teich}{{\mathfrak T}(M)}
\renewcommand{\o}{\operatorname}
\newcommand{\so}{{\o{SO}(2,1)^0}}
\newcommand{\SO}{{\o{SO}(2,1)}}
\newcommand{\rto}{\R^{2,1}}
\newcommand{\sgn}{{\o{sgn}}}
\newcommand{\vx}{{\mathsf x}}
\newcommand{\vv}{{\mathsf v}}
\newcommand{\vu}{{\mathsf u}}
\newcommand{\vw}{{\mathsf w}}
\newcommand{\xo}[1]{{\vx^0(#1)}}
\newcommand{\xp}[1]{{\vx^+(#1)}}
\newcommand{\xm}[1]{{\vx^-(#1)}}
\newcommand{\xpm}[1]{{\vx^{\pm}(#1)}}

\newcommand{\zo}{Z^1(G,\rto)}
\newcommand{\ho}{H^1(G,\rto)}
\newcommand{\hz}{H^1(\Z,\rto)}

\newcommand{\SOto}{\o{SO}(2,1)}
\newcommand{\Oto}{\o{O}(2,1)}
\newcommand{\SOoo}{\o{SO}^0(2,1)}
\newcommand{\isoo}{\o{Isom}^0(\rto)}
\newcommand{\iso}{\o{Isom}(\rto)}
\newcommand{\ee}{\mathcal{E}}
\newcommand{\tr}{\o{tr}}
\newcommand{\Hyp}{\o{Hyp}}
\newcommand{\sltr}{\o{SL}(2,\R)}
\newcommand{\lagltr}{\mathfrak{gl}(2,\R)}
\newcommand{\lasltr}{\mathfrak{sl}(2,\R)}
\newcommand{\laiso}{\mathfrak{o}(2,1)}
\newcommand{\la}{\langle}
\newcommand{\ra}{\rangle}
\newcommand{\Ad}{\o{Ad}}
\newcommand{\tit}{\tilde\iota_t}
\newcommand{\tg}{\tilde{g}}
\newcommand{\Hom}[1]{\o{Hom}(G,#1)}
\newcommand{\Homg}[1]{{\o{Hom}(G,#1)}/#1}
\newcommand{\Homl}{\Hom{\so}}
\newcommand{\Homa}{\Hom{\isoo}}
\newcommand{\Homs}{\Hom{\sltr}}
\newcommand{\Homsg}{\Homg{\sltr}}

\newtheorem{thm}{Theorem}
\newtheorem{lemma}[thm]{Lemma}
\newtheorem*{thm*}{Theorem}
\newtheorem*{remark*}{Remark}

\title[Flat Lorentz 3-Manifolds]
{Flat Lorentz 3-manifolds and cocompact Fuchsian groups}
\author[Goldman]{William M.\ Goldman}
\address{Department of Mathematics \\ 
University of Maryland\\ College Park, MD 20742 USA}
\email{wmg@math.umd.edu}
\author[Margulis]{Gregory A.\  Margulis}
\address{Department of Mathematics \\ 
10 Hillhouse Ave. \\
P.O. Box 208283 \\
Yale University\\ New Haven, CT 06520 USA}
\email{margulis@math.yale.edu}
\thanks{Both authors gratefully acknowledge partial support from NSF grants.}
\maketitle
% \begin{abstract}
% \tableofcontents
\section{Introduction}
Consider Minkowski $2+1$-space $\E$ and let $G\subset\so$ be a
discrete subgroup. Suppose that a group of affine
isometries of $\E$ with linear part $G$ acts properly and freely
on $\E$. In a remarkable preprint~\cite{Mess}, Geoffrey Mess proved
the following theorem:

\begin{thm*} 
$G$ is not cocompact in $\so$.
\end{thm*}
Mess deduces this result as part of a general theory of domains of
dependence in constant curvature Lorentzian 3-manifolds.  We give an
alternate proof, using an invariant introduced by
Margulis~\cite{Margulis1,Margulis2} and Teichm\"uller theory.
% (see also Drumm~\cite{Drumm1}).

We thank Scott Wolpert for helpful conversations concerning
Teichm\"uller theory. We also wish to thank Paul Igodt and the Algebra
Research Group at the Katholieke Universiteit Leuven at Kortrijk,
Belgium for their hospitality at the ``Workshop on Crystallographic
Groups and their Generalizations II'', where these results were
obtained.

\section{Background}
Let $\rto$ be a 3-dimensional real vector space with inner product
\begin{equation*}
\B(x,y) = x_1y_1 + x_2y_2 - x_3y_3.
\end{equation*}
The group of linear isometries of $\rto$ will be denoted by $\SOto$.
Let $\iso$ denote the group of {\em affine isometries, \/} that is, the
group of all transformations of the form
\begin{align*}
h: \rto &\longrightarrow \rto \\
x &\longmapsto g(x) + u
\end{align*}
where $g\in\Oto$ and $u\in\rto$.  We write $g=\L(h)$ and
$h=(g,u)$. Evidently $\iso$ is isomorphic to the semidirect product
$\Oto\ltimes\rto$ where $\rto$ denotes the vector group of {\em
translations\/} of $\E$.

Let $G\subset\Oto$ be a subgroup. An {\em affine deformation of $G$\/} is
a homomorphism $\phi:G\longrightarrow\iso$ such that $\L(\phi(g))=g$.
An affine deformation $\phi$ is {\em proper\/} if
the resulting action of $G$ by affine transformations on $\rto$ is a 
proper action. Write 
\begin{equation*}
\phi(g) = (g,u(g)).  
\end{equation*}
The condition that $\phi$ be a homomorphism is that the map
$u = u_\phi:G\longrightarrow\rto$ satisfy the 
{\em cocycle condition\/}
\begin{equation}\label{eq:coco}
u_\phi(g_1 g_2) = u_\phi(g_1) + g_1 u_\phi(g_2).
\end{equation}
A map $u:G\longrightarrow\rto$ satisfying \eqref{eq:coco} is called a
{\em cocycle\/} and the vector space of cocycles is denoted by $Z^1(G,\rto)$.

If $\phi_1,\phi_2$ are affine deformations of $G$ which are conjugate
by translation by $v\in\rto$, then the difference $u_{\phi_1} - u_{\phi_2}$
is the cocycle 
\begin{equation*}
\delta v: g \longmapsto v - g(v).
\end{equation*}
Such a cocycle is called a {\em coboundary.\/} The subspace of coboundaries
is denoted by $B^1(G,\rto)$. 
We say that $\phi_1,\phi_2$ are {\em translationally conjugate.\/}
Translational conjugacy classes of affine
deformations of $G$ correspond to elements in the {\em cohomology group\/}
\begin{equation*}
H^1(G,\rto) = Z^1(G,\rto)/B^1(G,\rto).
\end{equation*}

Suppose that $\phi:G\longrightarrow\iso$ is a proper affine deformation.
By Fried-Goldman~\cite{FriedGoldman}, the group $G$ is solvable 
or the linear part 
\begin{equation*}
\L\circ\phi: G \longrightarrow \Oto
\end{equation*}
is an isomorphism onto a discrete subgroup of $\Oto$. 
(Indeed, this conclusion is obtained 
for any proper affine action on $\R^3$.) 
The solvable groups are easily classified by embedding them as 
lattices in Lie  
subgroups which themselves act properly. 
When $G$ is not solvable, then interesting examples do exist
(Margulis~\cite{Margulis1,Margulis2}). Furthermore every {\em torsionfree
non-cocompact\/} discrete subgroup $G\subset\Oto$ for which
$H^1(G;\rto)\neq 0$   admits proper affine deformations (Drumm~\cite{Drumm2}).

Recall that an element of $\Oto$ is {\em hyperbolic\/} if it has three
distinct real eigenvalues.  A subgroup $G\subset\Oto$ is {\em purely
hyperbolic\/} if every element is hyperbolic. A cocompact discrete
subgroup contains a purely hyperbolic subgroup of finite index .

\section{An invariant of affine isometries}
In \cite{Margulis1,Margulis2}, Margulis defines an invariant
$\alpha_\phi: G \longrightarrow \R$
of an affine deformation $\phi$ of a purely hyperbolic subgroup
$G\subset\Oto$ as follows. We assume that $G\subset\so$. 
Choose a component $\Nn_+$ of the complement of $0$ in the 
lightcone. Since any element $g$ of $G$ is hyperbolic 
its three eigenvalues are distinct 
positive real numbers
\begin{equation*}
\lambda(g) < 1 < \lambda(g)^{-1}. 
\end{equation*}
\pagebreak

Choose an eigenvector $\xm{g}\in\Nn_+$ for $\lambda(g)$ and
an eigenvector $\xp{g}\in\Nn_+$ for $\lambda(g)^{-1}$, respectively.
Then there exists a unique eigenvector $\xo{g}$ for $g$ with eigenvalue $1$
such that:
\begin{itemize}
\item $\B(\xo{g},\xo{g}) = 1$;
\item $(\xm{g},\xp{g},\xo{g})$ is a positively oriented basis.
\end{itemize}
Notice that $\xo{g^{-1}}=-\xo{g}$.

If $\phi$ is an affine deformation corresponding to a cocycle $u$, 
then $\alpha_\phi$ is defined as:
\begin{align}\label{eq:defalpha}
\alpha_\phi: G & \longrightarrow \R \\
g &\longmapsto  \B(\xo{g},u(g)) \notag.
\end{align}
More generally, $\alpha_{\phi}(g) = \B(\xo{g},\phi(g)(x)-x))$ 
for any $x\in\E$. 
Furthermore $\alpha_{\phi}$ is a class function on $G$ and 
recently Drumm-Goldman~\cite{alpha} have proved that the mapping
\begin{align*}
H^1(G,\rto) & \longrightarrow \R^G \\ [u] & \longmapsto
\alpha_\phi
\end{align*}
is injective, that is, $\alpha$ is a complete invariant of the conjugacy class
of the affine deformation.

In \cite{Margulis1,Margulis2}, Margulis proved the following theorem
(see also Drumm~\cite{Drumm1}):

\begin{thm}[Margulis]\label{thm:SignCriterion}
Suppose that $G\subset\so$ is purely hyperbolic and let
$\phi:G\longrightarrow\iso$ be an affine deformation. If there exist
$g_1,g_2\in G$ such that $\alpha_\phi(g_1)>0 > \alpha_\phi(g_2)$,
then $\phi$ is not proper.
\end{thm}

Affine deformations defining free actions correspond to cocycles for
which $\alpha(g)\neq 0$ for $g\neq\Id$.  We shall say that a cocycle
$u$ is {\sl positive\/} (respectively {\sl negative\/}) if $\alpha(g) > 0$
(respectively $\alpha(g) < 0$) whenever $\Id\neq g\in G$.  Clearly $u$ is
positive if and only if $-u$ is negative.  We conjecture a converse to
Theorem~\ref{thm:SignCriterion}: {\em an affine deformation is proper
if and only if its cocycle is positive or negative.\/}

\section{Deformation-theoretic interpretation of $\alpha$}

We reduce the proof of Mess's theorem to facts about deformations of
hyperbolic Riemann surfaces.  Let $M$ be a surface with a complete
hyperbolic structure and $\pi=\pi_1(M)$ its fundamental group. A
representation $\phi:\pi\longrightarrow\so$ is {\em Fuchsian\/} if it
is an embedding onto a discrete subgroup of $\so$.  When $M$ is a
closed surface, the space of conjugacy classes of Fuchsian
representations $\phi:\pi\longrightarrow\so$ is an open subset of the
space of conjugacy classes of all representations, which identifies
with the {\em Teichm\"uller space\/} $\Teich$ of $M$. (See
Weil~\cite{Weil1,Weil2,Weil3}, \S VI of Raghunathan~\cite{Raghunathan}
for the general theory and Goldman~\cite{Goldman1,Goldman2} for the
case of surface groups.)  Its tangent space identifies with the
cohomology group $H^1(G,\rto)$ where $G=\phi(\pi)$.

Since the classical theory of Fuchsian groups is usually phrased in
terms of $\sltr$ (rather than $\SOto$), and since $2\times 2$ matrices
are more tractable than $3\times 3$ matrices, we work with $\sltr$.
The Lie groups $\sltr$ and $\SOto$ are {\em locally\/} isomorphic, but
not {\em globally\/} isomorphic. One model for the local isomorphism
is the adjoint representation, as follows.  The trace form of any
nontrivial representation (for example the Killing form) provides the
Lie algebra $\lasltr$ with a Lorentzian inner product invariant under
the adjoint representation. Thus $\lasltr$ is isometric to $\rto$; we
give an explicit orthogonal basis. In this way the adjoint representation
$\Ad:\sltr\longrightarrow\o{Isom}(\lasltr)$ defines a local
isomorphism $\rho:\sltr\longrightarrow\SOto$ of Lie groups.

The local isomorphism
$\rho:\sltr \longrightarrow \Oto $
is not injective --- its
kernel consists of the center $\{\pm\Id\}$ of $\sltr$.  
Nor is $\rho$ surjective --- 
its image is the identity component $\SOoo$ of $\Oto$.
Neither issue is problematic here, since 
purely hyperbolic discrete subgroups of 
$\SOto$ lift to subgroups of $\sltr$ (Abikoff~\cite{Abikoff},
Culler~\cite{Culler}, Kra~\cite{Kra}).
Let $G$ be a purely hyperbolic subgroup of $\SOto$, with inclusion
$\iota:G\hookrightarrow\SOto$. Then there exists a representation
$\tilde\iota:G\longrightarrow\sltr$ such that $\iota = \rho\circ\tilde\iota$.
Furthermore composition with the local isomorphism 
$\rho$ induces a covering space 
\begin{equation*}
\Homs\longrightarrow\Homa. 
\end{equation*}
Thus smooth paths in $\Homa$ lift to $\Homs$.  
Henceforth we suppress
$\tilde\iota$ (identifying $G$ with its image $\tilde\iota(G)$ in
$\sltr$) and consider paths in $\Homs$.

\section{$\lasltr$ and $\rto$}
For the calculations later, we now give a detailed description of the
local isomorphism $\rho$ derived from the adjoint representation.

For convenience, consider the Lie algebra $\lasltr$ with inner prouct
\begin{equation}\label{eq:defB}
\B(X,Y) := \frac12 \tr(XY).
\end{equation}
The basis
\begin{equation*}
e_1 = \bmatrix 1 & 0  \\ 0 & -1\endbmatrix,
e_2 = \bmatrix 0 & 1  \\ 1 & 0\endbmatrix,
e_3 = \bmatrix 0 & -1  \\ 1 & 0\endbmatrix.
\end{equation*}
is orthogonal with respect to $\B$ and satisfies
\begin{equation*}
\B(e_1,e_1) =  \B(e_2,e_2) = 1,\   \B(e_3,e_3) =  -1.
\end{equation*}
This provides an isometry of Lorentzian vector spaces
\begin{align*}
\psi: \lasltr &\longrightarrow  \rto \\
\bmatrix v_1 & v_2 \\ v_3 & -v_1 \endbmatrix
& \longmapsto 
\bmatrix v_1 \\ (v_2 + v_3)/2 \\ (-v_2 + v_3)/2 \endbmatrix. 
\end{align*}
With respect to this isometry the adjoint
representation defines a local isomorphism $\rho:\sltr \to \Oto$ satisfying:
\begin{equation*}
\psi(\Ad(g) v) = \rho(g)\psi(v) 
\end{equation*}
whenever $g\in\sltr$ and $v\in\lasltr$.
(In other words, $\psi:\lasltr_{\Ad}\to\rto$ 
is $\rho$-equivariant.)
Explicitly, 
\begin{align*}
\sltr &\stackrel{\rho}\to  \Oto \\
\bmatrix a & b \\ c & d \endbmatrix & \mapsto
%\bmatrix 1 + 2bc & -ab+ cd & ab + cd \\ 
\bmatrix 1 + 2bc & -ac+ bd & ac + bd \\ 
%-ac + bd & (a^2 - b^2 - c^2 + d^2)/2 &
- ab + cd & (a^2 - b^2 - c^2 + d^2)/2 & 
% (-a^2 + b^2 - c^2 + d^2)/2 \\ 
(  -a^2 - b^2 + c^2 + d^2)/2 \\
%ac + bd & (-a^2 - b^2 + c^2 + d^2)/2 &
 ab + cd & (-a^2 + b^2 - c^2 + d^2)/2 &
% (a^2 + b^2 + c^2 + d^2)/2   \endbmatrix
(  a^2 + b^2 + c^2 + d^2)/2 \endbmatrix
\end{align*}
(where $ad-bc=1$). Differentiation at $\Id\in\sltr$ (that is, at $a=d=1$,
$b=c=0$) gives the Lie algebra isomorphism
\begin{align*}
\lasltr & \longrightarrow  \laiso \\
\bmatrix v_1 & v_2 \\ v_3 & -v_1 \endbmatrix & \longmapsto \bmatrix 
0 & v_3-v_2 &  v_2 + v_3 \\ 
v_2 - v_3 & 0 & 2 v_1 \\ 
v_2 + v_3  & -2 v_1 & 0 \endbmatrix.
\end{align*}

An element $g\in\sltr$ is {\em hyperbolic\/} if it has two real distinct
eigenvalues, which are necessarily reciprocal. If $g$ has eigenvalues
$\mu,\mu^{-1}$ with $\vert\mu\vert < 1$, 
then $\rho(g)$ has eigenvalues $\lambda = \mu^2, 1,\mu^{-2}$. In particular
$g\in\sltr$ is hyperbolic if and only if $\rho(g)$ is hyperbolic.
There exists $f\in\sltr$ such that
\begin{equation*}
fg f^{-1} = g_0
\end{equation*}
where
\begin{equation*}
g_0 = \pm \bmatrix \mu & 0 \\ 0 & \mu^{-1} \endbmatrix  
\end{equation*}
and 
\begin{equation*}
0 < \mu < 1 < \mu^{-1}.
\end{equation*}
The eigenvectors of $g_0= \rho(g_0)$ are:
\begin{align*}
\xm{g_0} & = \psi\left( \bmatrix 0 & -2 \\ 0 & 0 \endbmatrix \right) 
= \bmatrix 0 \\ -1 \\ 1 \endbmatrix \\
\xp{g_0} & = \psi\left( \bmatrix 0 & 0 \\ 2 & 0 \endbmatrix \right) 
= \bmatrix 0 \\ 1 \\ 1 \endbmatrix \\
\xo{g_0} & = \psi\left( \bmatrix -1 & 0 \\ 0 & 1 \endbmatrix \right)  
= \bmatrix -1 \\ 0 \\ 0 \endbmatrix.
\end{align*}
The eigenvectors for $g$ are the images of the eigenvectors of $g_0$ under $f$.

Now we derive a formula for $\alpha(g)$ for an affine deformation $\phi$ which
is of the form $h = (\rho(g),\psi(v)(g))$ 
where $g\in G\subset\sltr$ and $v\in\lasltr$.
Suppose that $g\in\sltr$ is hyperbolic. We use the embedding
$\sltr\hookrightarrow\lagltr$.
Orthogonal projection 
\begin{align*}
\lagltr& \stackrel{\Pi}\longrightarrow\lasltr \\
g & \longmapsto  g - \frac{\tr(g)}2 \Id.
\end{align*}
maps $g_0$ to a diagonal matrix of trace zero. 
Dividing $\Pi(g_0)$ by 
\begin{equation*}
\sgn(\tr(g))\sqrt{-\det(\Pi(g_0))}  
\end{equation*}
gives the diagonal matrix
corresponding to $\xo{g_0}\in\rto$
(where  $\sgn(x)$ denotes the {\em sign\/} of a nonzero real number $x$).
Since $\tr(g_0) = \pm (\mu + \mu^{-1})$, 
\begin{equation*}
\det(\Pi(g_0)) = -(\mu-\mu^{-1})^2 = -\left(\tr(g_0)^2 -4\right)/4 
\end{equation*}
so
\begin{align*}
\sgn(\tr(g_0))
&\Pi(g_0) /\sqrt{-\det(\Pi(g_0))} \\  & = 
\sgn(\tr(g_0))
\left(g_0 - \frac{\tr(g_0)}2 \Id\right) \Big/
\left(\frac{\sqrt{\tr(g_0)^2-4}}{2}\right) \\  & = 
\bmatrix -1 & 0 \\ 0 & 1 \endbmatrix 
\end{align*}
corresponds to $\xo{g}$. Conjugation by $f$ gives the general formula
\begin{equation}\label{eq:xo}
\psi:
\sgn(\tr(g))\left(g - \frac{\tr(g)}2 \Id\right)
\Big/\left(\frac{\sqrt{\tr(g)^2-4}}2\right) 
\longmapsto \xo{g} 
% = \left(2 g - \tr(g) \Id\right)\big/\sqrt{\tr(g)^2-4}.
\end{equation}
From \eqref{eq:xo} follows a formula for $\alpha(g)$ in terms of traces.
Suppose that $G\subset\sltr$ is purely hyperbolic 
and $u\in Z^1(G,\lasltr)\cong\zo$. 
Taking the trace of the product of 
\eqref{eq:xo} with $u(g)$, and applying \eqref{eq:defalpha} and
\eqref{eq:defB} yields:
\begin{equation}\label{eq:alphat}
\alpha(g) =  
\sgn(\tr(g))
\frac{\tr\left(u(g)g\right)}{\sqrt{\tr(g)^2-4}}
\end{equation}

\section{Trace and displacement length}

Let $\Hyp$ denote the subset of $\sltr$ consisting of hyperbolic
elements.  The image of the trace function
$\tr:\Hyp\longrightarrow\R$ consists of the disjoint two intervals
$(-\infty,-2)$ and $(2,\infty)$.  Furthermore hyperbolic elements
$g\in\Hyp$ are determined up to conjugacy by their trace. 
In terms of hyperbolic geometry, $\tr(g)$ relates 
to the {\em displacement length $\ell(g)$,\/} 
that is, the minimum distance $g$ moves a point 
$x\in\rht$. This minimum is realized when $x$ lies in the
$g$-invariant geodesic, which is necessarily unique.  Equivalently
$\ell(g)$ is the length of the shortest homotopically nontrivial
closed curve in the quotient $\rht/\la g\ra$.  Such a shortest curve
is necessarily a simple closed geodesic.  Let $\tg\in\sltr$ be a lift of
$g\in\iso$ to $\sltr$, 
that is, $g=\rho(\tg)$.  Displacement length of $g$ relates
to $\tr(\tg)$ and the eigenvalue $0<\mu < 1$ by:
\begin{align*}
\ell(g) & = - 2 \log \mu \\
\vert\tr(\tg)\vert & =  2\cosh(\ell(g)/2) \\
\end{align*}
(the sign of $\tr(\tg)$ is ambiguous since $\ker(\rho)=\{\pm\Id\}$).
Since
\begin{equation}\label{eq:montrell}
\frac{d\vert\tr\vert}{d\ell} = \sinh(\ell/2) > 0
\end{equation}
trace depends monotonically on displacement length.

% The trace of the $3\times 3$ orthogonal
% matrix $\rho(g)$ relates to the $\tr(g)$ by:
% \begin{equation*}
% \tr(\rho(g)) = \tr(g)^2 - 1.
% \end{equation*}
% If $g$ has eigenvalues $\mu,\mu^{-1}$, then $\rho(g)$ has eigenvalues
% $\mu^2,1,\mu^{-2}$.

% The sign of $\alpha_\phi$ relates to trace as follows.
% Suppose that $\phi$ is an affine deformation of $G$ with cocycle $u$
% and let 
% \begin{equation*}
% \tit\in\Homs 
% \end{equation*}
% be a path starting at the inclusion $\iota:G\hookrightarrow\sltr$ 
% tangent to $u$. 
% Then  
% \begin{align*}
% \tr(\tit(g)) & = \tr\left(g \exp(t u(g) + O(t^2)\right) \\ & =
% \tr\left(g (\Id + t u(g) + O(t^2))\right) 
% \end{align*}
% so
% \begin{equation}\label{eq:derivtr}
% \frac{d}{dt}\Big|_{t=0}
% \tr(\tit(g)) = \tr\left(g u(g)\right) =
% \tr(\tit(g)) = \tr\left(\Pi(g) u(g)\right).
% \end{equation}

Associated to a cocycle $u\in\zo$
are real analytic paths $\tit$ in $\Homs$ of the form
\begin{equation*}
\tit(g) = g \exp\left(t u(g) + O(t^2)\right)
\end{equation*}
where $t$ is defined in an open interval $I_g$ containing zero.
(In general $I_g$ may depend on $g$.)
We say that the cocycle $u$ is {\em tangent\/} to the path $\tit$.

Given a path $\tit\in\Homs$ where $\tit(G)\subset\Hyp$, 
consider the two functions
\begin{align*}
\tau_g: I_g & \longrightarrow \R \\
t & \longmapsto \big|\tr\left(\tit(g)\right)\big|
\end{align*} and
\begin{align*}
L_g: I_g & \longrightarrow \R \\
t & \longmapsto \ell\left(\tit(g)\right).
\end{align*}
When $\tit$ corresponds to a path 
$\mu(t)$ in $\Teich$, then $L_g = \ell_g\circ\mu$ where $\ell_g:\Teich\to\R$
is the geodesic length function associated to $g$.

\begin{lemma}\label{lem:alphaLength}
Let $\phi$ be an affine deformation of $G$ corresponding to the cocycle
$u\in\zo$ and let $g\in G$. 
Suppose that $\mu(t)$ is a path in $\Teich$ tangent to $u$.
Then
\begin{equation}\label{eq:alphaLength}
\alpha_\phi(g) = L_g'(0).
\end{equation}
Furthermore $\alpha_\phi(g)$ and $\tau_g'(0)$ have the same sign.
\end{lemma}

% \begin{lemma}\label{lem:signs}
% Let $\phi$ be an affine deformation of $G$ corresponding to the cocycle
% $u\in\zo$ and let $g\in G$. 
% Suppose that $\mu(t)$ is a path in $\Teich$ tangent to $u$.
% The following conditions are equivalent:
% \begin{itemize}
% \item $\alpha_\phi(g)>0$;
% \item $\tau_g'(0) > 0$;
% \item $L_g'(0) > 0$.
% \end{itemize}
% \end{lemma}
\begin{proof} 
Let $\tit:G\longrightarrow\sltr$ be a smooth path of representations
starting at the inclusion $\iota$ corresponding to $\mu(t)$.
\begin{align*}
\tau_g'(0) & = \frac{d}{dt}\Big|_{t=0}\vert\tr\tit(g)\vert 		\\ 
& = \pm\frac{d}{dt}\Big|_{t=0}
				\tr \left(g (\exp(t u(g) + O(t^2)))\right) \\
& = \pm\frac{d}{dt}\Big|_{t=0}   
				\tr \left(g (\Id + t u(g) + O(t^2))\right) \\
& = \pm \tr\left(g u(g)\right) 				
\end{align*}
where the sign equals $\sgn(\tr(\tit(g))) = \sgn(\tr(\tilde\iota_0(g)))$. 
Applying \eqref{eq:alphat} to the last expression gives
\begin{equation}\label{eq:dtaualpha}
\tau_g'(0)  = \frac{\sqrt{\tr(g)^2 - 4}}2 \alpha(g). 
\end{equation}
Thus $\tau_g'(0)$ has the same sign as $\alpha(g)$ as claimed.

To prove \eqref{eq:alphaLength}, 
apply \eqref{eq:montrell} and the chain rule to obtain:
\begin{equation}\label{eq:dtauLength}
\tau_g'(0) = 
\sinh\left(\frac{L_g(0)}{2}\right) 
L_g'(0).
\end{equation}
Since 
\begin{equation*}
\sinh\left(\frac{L_g(0)}{2}\right) = \frac{\sqrt{\tr(g)^2 - 4}}2,
\end{equation*}
\eqref{eq:alphaLength} follows from
\eqref{eq:dtaualpha} and \eqref{eq:dtauLength}.
\end{proof}
Thus a cocycle is positive (respectively negative) in the sense of
Theorem~\ref{thm:SignCriterion} if and only if the corresponding
deformation in $\Teich$ increases (respectively decreases) lengths of
closed curves, to first order.

\section{Reduction to Teichm\"uller theory}

Suppose that $G\subset\sltr$ and $\phi:G\longrightarrow\iso$ is a
proper affine deformation.  By Theorem~\ref{thm:SignCriterion}, the
corresponding cocycle $u\in\zo$ is either positive or negative; by
replacing $u$ by $-u$ if necessary, we assume that $u$ is positive.

By Fried-Goldman~\cite{FriedGoldman}, $G$ is necessarily discrete and
is isomorphic to its image in the group of affine isometries.  Suppose
that $G$ is cocompact. By passing to a subgroup of finite index, we
may assume that $G$ is torsionfree. Then $G$ acts freely on 
the real hyperbolic plane $\rht$ and
since $G$ is discrete and cocompact, $\rht/G$ is a closed hyperbolic
surface $M$. Furthermore $G$ is isomorphic to the fundamental group
$\pi_1(M)$.  The representation $\tilde\iota$ corresponds to a point
$O$ in the Teichm\"uller space $\Teich$ and the cohomology class
$[u]\in\ho$ corresponds to a tangent vector $\upsilon$ to $\Teich$ 
at $O$.

\begin{lemma}\label{lem:monotone}
There exists a path $\mu(t)$ in $\Teich$, defined for all $0\le t<\infty$ 
starting at $O\in\Teich$ with tangent vector $\upsilon\in T_O\Teich$:
\begin{align}\label{eq:path}
\mu(0) & = O \\
\mu'(0) & = \upsilon\notag
\end{align}
such that, for each $g\in G$,  
the geodesic length function $\ell_g$ is convex
along $\mu(t)$.
% Furthermore
%\begin{equation*}
%\ell_g(\mu(t)) > \ell_g(\mu(0)).
%\end{equation*}
\end{lemma}

Assuming Lemma~\ref{lem:monotone} 
and that $u$ is positive, we obtain a contradiction.
Since $\alpha(g) > 0$, the directional derivative 
\begin{equation*}
\mu'(0)\ell_g = \upsilon\ell_g = L_g'(0)> 0  
\end{equation*}
by Lemma~\ref{lem:alphaLength}.
Convexity implies that $\mu'(t)\ell_g$ cannot decrease as 
$t\longrightarrow +\infty$. Thus 
\begin{equation*}
(\ell_g\circ\mu)'(t) = \mu'(t)\ell_g\ge \mu'(0)\ell_g = \alpha(g) > 0 
\end{equation*}
for all $t\ge 0$. In particular $\ell_g\circ\mu$ is monotone. Furthermore
\begin{equation}\label{eq:infinite}
\ell_g(\mu(t))\longrightarrow +\infty 
\text{~as~}t\longrightarrow +\infty,
\end{equation}
that is, {\em each closed geodesic on the hyperbolic
surface $\mu_t$ lengthens as $t\longrightarrow +\infty$.}

Such a path $\mu$ cannot exist for closed hyperbolic surfaces. 
Let $N>0$. Then for only finitely many conjugacy classes 
$F = \{[g_1],\dots,[g_m]\}$ 
in $G \cong \pi_1(M)$, 
the corresponding closed geodesics in $M$ have length $< N$. 
(Here $[g]$ denotes the conjugacy class of $g\in G$.)
For any $g\in G$ with $[g]\notin F$, the length function
$L_g(t) > L_g(0) \ge N$. Now consider $[g_i]\in F$. Let 
\begin{equation*}
\alpha_0 = \min_{1\le i\le m}\alpha(g_i)> 0.  
\end{equation*}
Convexity, together with \eqref{eq:alphaLength} implies that 
\begin{equation*}
L_{g_i}(t) \ge  L_{g_i}(0) + t \alpha(g_i) \ge t\alpha_0.
\end{equation*}
Hence, for $t > N/\alpha_0$, 
\begin{equation*}
L_g(t) = \ell_g(\mu_t) > N 
\end{equation*}
for all $g\in G-\{\Id\}$.

However, for any closed hyperbolic surface $M$ there exists
a simple closed geodesic of length at most $2\log(2-2\chi(M))$
(Lemma~5.2.1 of Buser~\cite{Buser}). 
Taking $N > 2\log(2-2\chi(M))$, we obtain the desired contradiction. \qed
% (This lemma follows from the fact  --- for example,
% Lemma~4.1.5 of Buser~\cite{Buser} --- that the injectivity radius equals
% half the length of the shortest closed geodesic.)  

\begin{proof}[Proof of Lemma~\ref{lem:monotone}]
Here are two constructions for $\mu$, the first based on the
Riemannian geometry of $\Teich$ with the Weil-Petersson metric and the
second based on Thurston's earthquake flows.

Let $\mu(t)$ be the Weil-Petersson geodesic satisfying \eqref{eq:path}.  
By Corollary~4.7 of Wolpert~\cite{Wolpert2}, 
the geodesic length function $\ell_g$ is strictly
convex along $\mu(t)$ and the directional derivative 
$\upsilon\ell_g > 0$, for any $g\in G -\{1\}$.  
Therefore $\ell_g\circ\mu(t)$ is monotonically increasing for $t>0$.

However, in general the Weil-Petersson metric is geodesically
incomplete (Chu~\cite{Chu}, Wolpert~\cite{Wolpert3}), so that $\mu(t)$
is only defined for $t_1<t<t_2$ where $t_1< 0 < t_2$.  We show this is
impossible under our assumptions on $\mu'(0)=\upsilon$.

By Mumford's compactness theorem (Mumford~\cite{Mumford},
Harvey~\cite{Harvey}, 2.5.1 or Buser~\cite{Buser}, 6.6.5), 
the subspace of moduli space consisting of hyperbolic surfaces
whose injectivity radius is larger than any positive constant is compact. 
An incomplete geodesic on a Riemannian manifold
must leave every compact set. Therefore, if the
Weil-Petersson geodesic $\mu(t)$ cannot be extended to $t_2 <\infty$,
then
\begin{equation*}\label{eq:nonsense}
\lim_{t\to t_2}\ \inf_{g\in G-\{\Id\}}\ell_g(\mu(t)) = 0,
\end{equation*}
contradicting monotonicity of $\ell_g$.

Hence $\mu(t)$ is defined for all $t<\infty$.
As above, convexity implies \eqref{eq:infinite}.

Alternatively, take $\mu$ to be the earthquake path
% $\Ee_\lambda$ 
introduced by Thurston (see
Kerckhoff~\cite{Kerckhoff1,Kerckhoff2} and 
Thurston~\cite{Thurston1}).
For the given tangent vector $\upsilon$, there exists a unique
measured geodesic lamination $\lambda$ such that the corresponding
earthquake path $\mu(t) = \ee_\lambda(t)$ satisfies
\eqref{eq:path}  (Kerckhoff~\cite{Kerckhoff2}, Proposition~2.6).
By Kerckhoff~\cite{Kerckhoff1} (see also Wolpert~\cite{Wolpert1}),
each length function $\ell_g$ is convex along the earthquake
path $\ee_\lambda$, implying \eqref{eq:infinite}. Indeed,
$\ell_g$ is strictly convex along $\mu$
since the lamination $\lambda$ {\em fills up\/} $M$ --- that is,
every nonperipheral simple closed curve $\sigma$ intersects $\lambda$. For 
otherwise $\ell_\sigma$ would be constant along $\mu$, contradicting
\begin{equation*}
\frac{d}{dt}\Big|_{t=0} \ell_\sigma\circ\mu(t)  > 0. 
\end{equation*}
% Corollary 4.7 of \cite{Wolpert2}, a geodesic length function is strictly
% convex along a Weil-Petersson geodesic.
\end{proof}
\begin{remark*}\rm
Another proof, closer in spirit to the proof in \cite{Mess}, involves the
density of simple closed curves in the projective measured lamination
space. Let $\Ss$ denote the set of isotopy classes of simple closed
curves on $M$ and let $\PL$ denote Thurston's space of projective
equivalence classes of measured geodesic laminations on $M$. Since 
\begin{align*}
\ML & \longrightarrow T_O\Teich \\
\lambda &\longmapsto \Ee_{\lambda}'(0)
\end{align*}
is a homeomorphism (Proposition~2.6 of \cite{Kerckhoff2}), there
exist $\lambda\in\ML$ satisfying $\Ee_\lambda'(0)=\upsilon\neq 0$.
Theorem~5.1 of \cite{Thurston2} implies
\begin{align*}
\PL &\longrightarrow  T_O^*\Teich \\
[\lambda] & \longmapsto d\log\ell_\lambda
\end{align*}
is an embedding onto a convex sphere in $T_O^*\Teich$
(where $\ell_\lambda(N)$ denotes the {\em length\/} of the lamination
$\lambda$ as measured in $N$).
Since $\Ss$ is dense in $\PL$,
there exist $\gamma_1,\gamma_2\in\Ss$ such that
\begin{align*}
(d\log\ell_{\gamma_1})(\lambda) &> 0 \\
(d\log\ell_{\gamma_2})(\lambda) & < 0. 
\end{align*}
Let $g_1,g_2\in\pi_1(M)$ correspond to $\gamma_1,\gamma_2$ respectively.
Then 
\begin{equation*}
\upsilon L_{g_1} >0, \upsilon L_{g_2} <0, 
\end{equation*}
contradicting Theorem~\ref{thm:SignCriterion} and Lemma~\ref{lem:alphaLength}.
\end{remark*}
\begin{remark*}\rm
Mess's original proof uses Lorentzian geometry, and in particular the
theory of domains of dependence in constant curvature Lorentzian space
forms developed in \cite{Mess} and Scannell~\cite{Scannell}.  As part
of his general theory, Mess shows that any affine deformation
sufficiently near the holonomy of a complete flat Lorentz 3-manifold
is the holonomy of a complete flat Lorentz 3-manifold, that is, the
nearby action is also proper and free.
The cocycle $u$ corresponds to the velocity vector to an earthquake path
$\ee_\lambda$ along a measured geodesic lamination $\lambda$, and 
$\lambda$ is approximated by a
{\em finite measured geodesic lamination,\/} that is, a disjoint union
of simple closed geodesics.
However for a finite lamination, the corresponding group action
is not free (elements of $G$ corresponding to curves disjoint from $\lambda$
have fixed points), a contradiction.
\end{remark*}

\makeatletter \renewcommand{\@biblabel}[1]{\hfill#1.}\makeatother

\end{document}